\documentclass[12pt,leqno]{article}
\usepackage{amsmath, amsthm, amsfonts, amssymb, color}
\usepackage{mathrsfs}
\theoremstyle{definition}

\newcommand{\scr}[1]{\mathscr #1}
\definecolor{wco}{rgb}{0.5,0.2,0.3}

\numberwithin{equation}{section} \theoremstyle{remark}

\newcommand{\ua}{\uparrow}

\title{{\bf Modified Curvatures on Manifolds with Boundary and Applications}\footnote{Supported in
 part by WIMCS and SRFDP.}
}
\author{
{\bf Feng-Yu Wang}\\
\footnotesize{School of Mathematical Sciences,
Beijing Normal
University, Beijing 100875, China}\\
\footnotesize{and}\\ \footnotesize{Department of Mathematics,
Swansea University, Singleton Park, SA2 8PP, UK}\\ \footnotesize{Email: wangfy@bnu.edu.cn;
F.Y.Wang@swansea.ac.uk}}
\begin{document}
\def\R{\mathbb R}  \def\ff{\frac} \def\ss{\sqrt} \def\BB{\mathbb
B}
\def\N{\mathbb N} \def\kk{\kappa} \def\m{{\bf m}}
\def\dd{\delta} \def\DD{\Delta} \def\vv{\varepsilon} \def\rr{\rho}
\def\<{\langle} \def\>{\rangle} \def\GG{\Gamma} \def\gg{\gamma}
  \def\nn{\nabla} \def\pp{\partial} \def\tt{\tilde}
\def\d{\text{\rm{d}}} \def\bb{\beta} \def\aa{\alpha} \def\D{\scr D}
\def\E{\mathbb E} \def\si{\sigma} \def\ess{\text{\rm{ess}}}
\def\beg{\begin} \def\beq{\begin{equation}}  \def\F{\scr F}
\def\Ric{\text{\rm{Ric}}} \def\Hess{\text{\rm{Hess}}}\def\B{\scr B}
\def\e{\text{\rm{e}}} \def\ua{\underline a} \def\OO{\Omega} \def\sE{\scr E}
\def\oo{\omega}     \def\tt{\tilde} \def\Ric{\text{\rm{Ric}}}
\def\cut{\text{\rm{cut}}} \def\P{\mathbb P} \def\ifn{I_n(f^{\bigotimes n})}
\def\C{\scr C}      \def\aaa{\mathbf{r}}     \def\r{r}
\def\gap{\text{\rm{gap}}} \def\prr{\pi_{{\bf m},\varrho}}  \def\r{\mathbf r}
\def\Z{\mathbb Z} \def\vrr{\varrho} \def\ll{\lambda}
\def\L{\scr L}\def\Tt{\tt} \def\TT{\tt}\def\II{\mathbb I}
\def\i{{\rm i}}\def\Sect{{\rm Sect}}

\maketitle
\begin{abstract}  To study  the reflecting diffusion processes on manifolds with boundary, 
some new curvature operators are introduced by using the Bakry-Emery curvature and the second fundamental form. As applications,   
 the gradient estimates, log-Harnack inequality and Poincar\'e/log-Sobolev inequalities are investigated for the Neumann semigroup on manifolds with boundary.
\end{abstract} \noindent
 \vskip 2cm

\section{Introduction}
The Bakry-Emery curvature condition \cite{BE} has played a crucial role in the study of diffusion semigroups on Riemannian manifolds. When the reflecting diffusion processes are considered on a manifold with boundary, both the curvature of the generator and the second fundamental form of the boundary have to be taken into account, see \cite{Hsu,W05, W07, W10, W10b} and references within.  It has been observed in \cite{W09b,W09c} that the curvature and the second fundamental form  play  essentially different roles in the study of functional inequalities for the reflecting diffusion processes, so that they can not be compensated each other. Moreover, since the geometry of the boundary works to a reflecting diffusion process only when the process reaches the boundary, the second fundamental form appears in the study  as   integrals w.r.t. the local time of the process on the boundary  (see \cite{Hsu, W10, W10b}). To avoid using the local time which is in general  less explicit,  we aim to derive explicit results for the reflecting diffusion processes by using modified curvature tensors consisting of   the Bakry-Emery curvature and information from the boundary.

Let $(M,\<\cdot,\cdot\>)$ be a $d$-dimensional complete connected Riemannian manifold with boundary $\pp M$. Let $P_t$ be the  semigroup of the reflecting diffusion  generated by $L:=\DD+Z$ for some $C^1$-vector field $Z$ on $M$. We assume that the reflecting diffusion process generated by $L$ is non-explosive, so that $P_t$ is a Markov semigroup. According to e.g. Lemma \ref{L3.0} below, $P_t$ is the Neumann semigroup generated by $L$; that is, for a reasonable reference function $f$, the times-space function $u:= P_\cdot f$ solves the Neumann problem
$$ \pp_t u= Lu,\ \ Nu|_{\pp M}:=\<N,\nn u\>|_{\pp M}=0, u(0,\cdot)=f,$$ where $N$ is the inward unit normal vector field of $\pp M$.

Recall that for any $f\in C^\infty(M)$, the Bochner-Weitzenb\"ock formula implies

\beq\label{1.1} \GG_2(f):= \ff 1 2 L|\nn f|^2 -\<\nn Lf, \nn f\>=(\Ric-\nn Z)(\nn f,\nn f) +\|\Hess_f\|_{HS}^2,\end{equation} where $\Ric$ and $\Hess$ stand for the Ricci curvature and the Hessian tensor respectively, and $\|\cdot\|_{HS}$ is the Hilbert-Schmidt norm. Consequently, for a function $g$ on $M$, the Bakry-Emery curvature condition
$$\GG_2(f)\ge g |\nn f|^2, f\in C^\infty(M)$$ is equivalent to $\Ric-\nn Z\ge g.$ Here and in the sequel, for a 2-tensor $\mathbf T$ and a function $g$, $\mathbf T \ge g$ means that
$\mathbf T(X,X)\ge g |X|^2$ holds for any  $X\in TM$, the tangent space of $M$.

Next,  The second fundamental form is a two-tensor on $\pp M$ given by
$$\II(X,Y):= -\<\nn_X N, Y\>,\ \ X,Y\in T\pp M,$$ where $T\pp M$ is the tangent space of $\pp M$.

 Now, for any strictly positive $\phi\in C^2(M)$, we introduce a family of modified curvature tensors
$$\Ric_Z^{\phi,p}:= \Ric -\nn Z -\ff 1 p (\phi^p L \phi^{-p})\<\cdot,\cdot\>,\ \ p>0.$$ To ensure that these tensors  contain also information from the boundary, the function $\phi$ will be taken from the class
$$\D:= \big\{\phi\in C_b^2(M): \inf\phi=  1, N\phi=0,  \II\ge - N\log \phi\big\}.$$  Note that  for a vector $X$ and a function $f$ we write  $Xf=\<X,\nn f\>$, and  conditions on $N$ and  $\II$ are automatically restricted to  $\pp M$ and $T\pp M$. If $\II$ and the sectional curvatures of $M$ are bounded and the injectivity radius of the boundary is positive, then the class $\D$ is non-empty, see \cite[Page 1436]{W07} for construction of $\phi$ using the distance function to the boundary. We also remark that the condition $\inf\phi=1$ in the definition of class $\D$ is not essential but for convenience, since our results (see Theorem \ref{T1.1} below) do not change if one replaces $\phi$ by $c\phi$ for a   constant $c>0.$

To construct the reflecting diffusion process,
let $B_t$ be the $d$-dimensional Brownian motion on a complete filtered probability space $(\OO, \{\F_t\}_{t\ge 0},\P)$. Then the reflecting diffusion process $X_t$ and its local time $l_t$ on $\pp M$ can be constructed by solving the Stratonovich stochastic differential equation on $M$:
\beq\label{E} \d X_t = \ss 2\,u_t\circ\d B_t +Z(X_t)\d t + N(X_t)\d l_t,\end{equation} where $u_t$ is the horizontal lift of $X_t$ on the frame bundle $O(M)$; i.e. $u_t$ satisfies
$$\d u_t= \mathbf H_{u_t}\circ \d X_t$$ for a fixed initial data $u_0\in O_{X_0}(M)$ and the horizontal lift  $\mathbf H$ from $TM$  to $T OM$ (the tangent space of $O(M)$). We have 
$$P_t f(x)= \E^x f(X_t),\ \ t\ge 0, x\in M, f\in \B_b(M),$$ where $\B_b(M)$ is the set of all bounded measurable functions on $M$, and $\E^x$ is the expectation taken  for the  process $X_t$ starting at point $x$.

 Let $X_t^\phi$ be the reflecting diffusion process generated by $$L^\phi:= L-2\nn \log\phi.$$
Since $X_t$ is non-explosive, so is $X_t^\phi$ provided $\nn\log \phi$ is bounded.  Below is the main result of the paper, which provides sharp gradient estimates of $P_t$ without using the local time. Let   $\rr$ be the Riemannian distance on $M$, i.e. for any $x,y\in M$, $\rr(x,y)$ is the length of the shortest curve on $M$ which links $x$ and $y$. For a fixed point $o\in M$, let $\rr_o=\rr(o,\cdot)$. We will need the following technical assumption:
\paragraph{(A)} \emph{ Either $|\nn\phi|\cdot |Z|$ is bounded, or $\Ric$ is bounded below and $|Z|\le \psi\circ\rr_o$ holds for some strictly positive function $\psi\in C([0,\infty))$ with $\int_0^\infty \ff 1 {\psi(r)}\d r =\infty.$}

\beg{thm}\label{T1.1} Let $\phi\in \D$ such that {\bf (A)} holds. Then for any $K\in C_b(M)$,  the following statements are equivalent to each other:
\beg{enumerate} \item[$(1)$] $\Ric_Z^{\phi,1}\ge K$; \item[$(2)$] For any $f\in C_b^1(M)$,
$$|\nn P_t f(x)| \le \ff 1 {\phi(x)}\E^x\Big\{ (\phi |\nn f|)(X_t) \e^{-\ss 2 \int_0^t \<u_s^{-1}\nn\log \phi (X_s),\, \d B_s\> -\int_0^t (K+|\nn\log \phi|^2)(X_s)\d s}\Big\}$$ holds for $ t\ge 0$ and $ x\in M;$ \item[$(3)$] For any $f\in C_b^1(M)$
 and $t\ge 0,\
 |\nn P_t f(x)| \le \dfrac 1 {\phi(x)}\E^x\Big\{ (\phi |\nn f|)(X_t^\phi) \e^{-\int_0^t  K (X_s^\phi)\d s}\Big\}.$  \end{enumerate}\end{thm}

As applications of Theorem \ref{T1.1}, we have the following explicit gradient/Poincar\'e/Harnack type inequalities for $P_t$.

\beg{cor}\label{C1.2} Let $\phi\in\D$ such that {\bf (A)} and
$ \Ric_Z^{\phi,2}\ge K_\phi$  hold  for some constant $K_\phi$. Then:
\beg{enumerate} \item[$(1)$] $\phi^2 |\nn P_t f|^2 \le  \e^{-2K_\phi t} P_t (\phi |\nn f|)^2$ holds for   any $f\in C_b^1(M)$ and $t\ge 0.$
\item[$(2)$] For any measurable function $f\ge 1$, the log-Harnack inequality
$$P_t\log f(y)\le \log P_t f(x) +\ff{\|\phi\|_\infty^2K_\phi \rr(x,y)^2}{2(\e^{2K_\phi t}-1)},\ \ t\ge 0, x,y\in M$$ holds.   \item[$(3)$] $P_t f^2\le (P_t f)^2 + \dfrac{\|\phi\|_\infty^2(1-\e^{-2K_\phi t})}{K_\phi} P_t |\nn f|^2 $   holds for   any $f\in C_b^1(M)$ and $t\ge 0.$
\item[$(4)$] $P_t f^2\ge (P_t f)^2 + \dfrac{\e^{2K_\phi t}-1}{\|\phi\|_\infty^2 K_\phi}  |\nn P_t f|^2 $   holds for   any $f\in C_b^1(M)$ and $t\ge 0.$
\end{enumerate} \end{cor}

\paragraph{Remark.} (a) The log-Harnack inequality was introduced in \cite{RW} for diffusion semigroups on Hilbert spaces with non-constant diffusion coefficients, which implies heat kernel bounds and the HWI (energy/cost/information) inequality. This inequality has been established in \cite[Section 5]{W10} on manifolds with boundary by using exponential estimates on the local time.

(b) Let $Z=\nn V$ for some $V\in C^2(M)$ such that $\mu(\d x):=\e^{V(x)}\d x$ is a probability measure, where $\d x$ stands for the volume measure on $M$. Then $P_t$ is symmetric in $L^2(\mu)$ and $P_t f\to \mu(f)$ in $L^2(\mu)$ as $t\to\infty$ for any $f\in L^\infty(M),$ where $\mu(f):=\int_M f\d\mu$. If $K_\phi>0$, by letting $t\to \infty$ in Corollary \ref{C1.2}(3) we obtain the Poincar\'e inequality

\beq\label{P} \mu(f^2)\le \mu(f)^2 + \ff{\|\phi\|_\infty^2}{K_\phi}\mu(|\nn f|^2),\ \ f\in C_b^1(M).\end{equation}  

\

Below we establish the corresponding log-Sobolev inequality and the HWI inequality, which generalize the  existing ones in the case without boundary. In particular, if $\pp M$ is convex  we may take $\phi\equiv 1$ so that Corollary \ref{C1.3}(1) goes back  to the Bakry-Emery criterion while Corollary \ref{C1.3}(2) reduces to the HWI inequality derived in \cite{OV} and \cite{BGL} on manifolds without boundary. Moreover, Corollary \ref{C1.3}(3) provides explicit heat kernel bounds.

\beg{cor}\label{C1.3} Let $Z=\nn V$ for some $V\in C^2(M)$ such that $\mu(\d x):=\e^{V(x)}\d x$ is a probability measure. Let $\phi\in \D$  such that {\bf (A)} and
$\Ric_Z^{\phi,2}\ge K_\phi$  hold for some constant $K_\phi$. \beg{enumerate} \item[$(1)$] If $K_\phi>0$  then
$$\mu(f^2\log f^2)\le \mu(f^2)\log \mu(f^2) + \ff{2\|\phi\|_\infty^6}{K_\phi}\mu(|\nn f|^2),\ \ f\in C_b^1(M), \mu(f^2)=1.$$
\item[$(2)$] If $K_\phi\le 0$ then
$$\mu(f^2\log f^2)\le 2\|\phi\|_\infty^4   \ss{\mu(|\nn f|^2)} W_2^\rr(f^2\mu, \mu)
 -\ff {\|\phi\|_\infty^2K_\phi}2 W_2^\rr(f^2\mu,\mu)^2$$ holds for any $f\in C_b^1(M)$ with $\mu(f^2)=1.$
 \item[$(3)$] Let $p_t(x,y)$ be the heat kernel of $P_t$ w.r.t. $\mu$. Then
 \beg{equation*}\beg{split} & p_t (x,y)\ge \exp\bigg[-\ff{\|\phi\|_\infty^2 K_\phi \rr(x,y)^2}{2(\e^{K_\phi t}-1)}\bigg],\ \ \text{and}\\
 & \int_M p_t(x,z)\log\ff{p_t(x,z)}{p_t(y,z)}\,\mu(\d z)\le \ff{\|\phi\|_\infty^2K_\phi\rr(x,y)^2}{2(\e^{2K_\phi t}-1)}\end{split}\end{equation*} hold for all $t>0$ and $x,y\in M.$\end{enumerate} \end{cor}

We will prove Theorem \ref{T1.1} and its Corollaries in Sections 2 and 3 respectively. Indeed, Section 2 proves  more than Theorem \ref{T1.1}: a result more general than the equivalence of (2) and (3) in Theorem \ref{T1.1} is proved (see Proposition \ref{P2.3} below), and   $\II\ge -N\log \phi$ is  deduced  from Theorem  \ref{T1.1}(3) for a class of manifolds with boundary including compact ones (see Proposition \ref{P2.4} below).

\section{Proof of Theorem \ref{T1.1}}
\subsection{From (1) to (2)}  We will make use of \cite[Proposition A.2]{W10}, for which we have to confirm that for any $f\in C_b^1(M)$ and $T>0$, $P_\cdot f$ is bounded on $[0,T]\times M$.  To this end, we  first   extend a result in \cite{W07} to make the boundary convex by using conformal changes of metric, then prove the boundedness of the gradient by following the line of \cite{W10}.

\beg{lem}[\cite{W07}]\label{L2.1} Let $\phi\in C^2(M)$ be strictly positive with $N\phi =0.$ If $\II\ge -N\log\phi$ then $\pp M$ is convex under the metric $\<\cdot,\cdot\>':=\phi^{-2}\<\cdot,\cdot\>.$\end{lem}
\beg{proof} Let $\nn'$ be the Levi-Civita connection for the metric $\<\cdot,\cdot\>'.$  We have (see \cite[Theorem 1.159(a)]{Bess})
$$\nn'_XY= \nn_XY -\<X,\nn\log \phi\>-\<Y,\nn\log \phi\>X +\<X,Y\>\nn\log\phi.$$
Since $\<X,\nn\phi\>=0$ for $X\in TM$ and noting that the inward unit normal vector field of $\pp M$ under the metric $\<\cdot,\cdot\>'$ is $N':=\phi N$, we obtain
\beg{equation*}\beg{split} -\<\nn'_XN', X\>' &= \phi^{-2} \<N',\nn\log \phi\>|X|^2-\phi^{-2}\<\nn_XN', X\> \\
&= \phi^{-1}\big(\II(X,X) +(N\log \phi)|X|^2\big)\ge 0.\end{split}\end{equation*}\end{proof} The second lemma is essentially due to \cite{W10}.  But we are using a different condition.

\beg{lem} \label{L2.2} Let $\phi\in\D$ such that $\Ric_Z^{\phi,1}$ is bounded below and {\bf (A)} holds,  then
$$\|\nn P_tf\|_\infty\le c\e^{ct} \|\nn f\|_\infty,\ \ t\ge 0, f\in C_b^1(M)$$ holds for some constant $c>0$. \end{lem}

\beg{proof} (a) Let  $|\nn \phi|\cdot|Z|$ be  bounded. By Lemma \ref{L2.1}, $\pp M$ is convex under the metric $\<\cdot,\cdot\>':=\phi^{-2}\<\cdot,\cdot\>.$
Let $\DD'$ and $\Ric'$ be the Laplacian and the Ricci curvature for the metric $\<\cdot,\cdot\>'.$ We have
$$\phi^2L= \DD' + (d-2) \phi\nn\phi +\phi^2 Z=: \DD' +Z'.$$ Since $|\nn\log\phi|$  and $|Z|\cdot|\nn\log\phi|$ are bounded, according to the calculations in the proof of \cite[Lemma A.4]{W10}, $\Ric'-\nn'Z'$ is bounded below and
the desired gradient inequality holds for some constant $c$ and all $f\in C_0^1(M)$. By an approximation argument as in the proof of \cite[Proposition A.2]{W10}, this inequality holds for all $f\in C_b^1(M)$.

(b) Let $\Ric$ be bounded below and $|Z|\le  \psi\circ\rr_o$ holds for some strictly positive function $\psi\in C([0,\infty))$ with $\int_0^\infty \ff 1 {\psi(r)}\d r =\infty.$ Then
$$\varphi:= \int_0^{\rr_o}\ff{\d r}{\psi(r)}$$ is a compact function, i.e. $\{\varphi\le r\}$ is compact for any constant $r$. For any $n\ge 1$, let $h_n= (2-\varphi/n)^+\land 1$ and $Z_n= h_n Z$. Then $|\nn\phi|\cdot |Z_n|$ is bounded and by our condition
$$\Ric_{Z_n}^{\phi,1}= h_n \Ric_Z^{\phi,1} + (1-h_n)\{\Ric- \phi\DD \phi^{-1}\<\cdot,\cdot\>\} +(\nn h_n)\otimes Z \ge K_0$$ holds for some constant $K_0$ independent of $n$. Therefore, by (a) there exists a constant $c_0>0$ such that
$$\|\nn P_t^{(n)} f\|_\infty\le c_0\e^{c_0 t} \|\nn f\|_\infty,\ \ t\ge 0, n\ge 1, f\in C_b^1(M),$$ where $P_t^{(n)}$ is the  semigroup of the reflecting diffusion process generated by $\DD +Z_n$. Since $P_t^{(n)}f\to P_t f$ as $n\to\infty$, we obtain
$$\ff{|P_t f(x)-P_t f(y)|}{\rr(x,y)}=\ff{|P_t^{(n)} f(x)-P_t^{(n)} f(y)|}{\rr(x,y)}\le c_0\e^{c_0t} \|\nn f\|_\infty,\ \ x,y\in M.$$ Therefore,  the desired gradient inequality holds. \end{proof}

Finally, the following lemma is an extension of    \cite[Lemma 2.3]{ATW} where the initial points is outside the boundary.

\beg{lem}\label{L2.3} For any $x\in M$ and $r_0>0$, there exists a constant $c>0$ such that
$$\P(\si_r\le s)\le \e^{-c r^2/t},\ \ r\in [0,r_0], t>0$$ holds, where $\si_r=\inf\{s\ge 0: \rr(X_s, x)\ge r\}$ and $X_s$ is the reflecting diffusion process generated by $L$ with $X_0=x$.\end{lem}

\beg{proof} Let $\phi\in C_b^\infty(M)$ such that $\phi\ge 1$ and $\pp M$ is convex in $B(x,r_0)$ under the metric $\<\cdot,\cdot\>':=\phi^{-2}\<\cdot,\cdot\>.$ As in the proof of Lemma \ref{L2.2}, we have $L= \phi^{-2}(\DD'+Z')$. Let $\rr'$ be the Riemannian distance function to $x$ induced by the metric $\<\cdot,\cdot\>'$. By taking smaller $r_0$ we may and do assume that $(\rr')^2\in C^\infty(B(x, 2 r_0))$. By the convexity of the boundary under the new metric and using the It\^o formula, we obtain
$$\d \rr'(X_t)^2 \le \ss 2 \phi^{-1} (X_t) \rr'(X_t)\d b_t +c_1\d t,\ \ t\le\si_{r_0}$$ for some constant $c_1>0$ and an one-dimensional Brownian motion $b_t$. Due to this inequality, the remainder of the proof is completely similar to that of \cite[Lemma 2.3]{ATW}.\end{proof}

\beg{proof}[Proof of (1) implying (2)]  Let $\phi\in \D$ such that (1) holds.  Since (1)  implies $\Ric-\nn Z\ge K+ \phi L\phi^{-1}$ while $\phi\in\D$ ensures $\II\ge -N\log\phi$, according to \cite[Proposition A.2]{W10} we have
$$ |\nn P_t f(x)| \le \E^x \Big\{|\nn f|(X_t) \e^{-\int_0^t (K+\phi L\phi^{-1})(X_s)\d s +\int_0^t  N\log \phi(X_s)\d l_s}\Big\} $$ provided $\E \e^{\ll l_t} <\infty$ holds for any $\ll, t>0.$ In general, the proof of \cite[Proposition A.2]{W10} implies that
$$\eta_s:= |\nn P_{t-s} f| (X_s) \e^{-\int_0^s (K+\phi L\phi^{-1})(X_r) \d r +\int_0^s N \log \phi(X_r) \d l_r},\ \ s\in [0,t]$$ is a local submartingale; that is, letting
$$\tau_n=\inf\{s \ge 0: l_s\lor \rr_o(X_s)\ge n\}$$ which goes to $\infty$ as $n\to\infty$, $\{\eta_{s\land\tau_n}\}_{s\in [0,t]}$ is a submartingale for each $n\ge 1$. So,
\beq\label{a1} |\nn P_t f(x)| \le \E^x \Big\{|\nn P_{(t-\tau_n)^+}f|(X_{t\land \tau_n}) \e^{-\int_0^{t\land\tau_n}  (K+\phi L\phi^{-1})(X_s)\d s +\int_0^{t\land\tau_n}  N\log \phi(X_s)\d l_s}\Big\} \end{equation} holds for $n\ge 1$.

On the other hand, by (\ref{E}) and the It\^o formula, we have
$$\d \log \phi(X_s) = \ss 2  \<u_s^{-1} \nn\log\phi(X_s), \d B_s\> + L\log \phi(X_s)\d s + N\log\phi(X_s)\d l_s.$$ Then
$$\int_0^{t\land\tau_n} N\log \phi(X_s)\d l_s = \log\ff{\phi(X_{t\land\tau_n})}{\phi(x)} -\ss 2 \int_0^{t\land\tau_n} \<u_s^{-1} \nn\log\phi(X_s), \d B_s\> -\int_0^{t\land\tau_n} L\phi(X_s)\d s.$$ Combining this with (\ref{a1}) and noting that
$$ \phi L \phi^{-1} - L\log\phi = -|\nn \log\phi|^2,$$ we obtain
\beg{equation*}\beg{split} |\nn P_t f(x)| \le \ff 1 {\phi(x)}\E^x\Big\{& (\phi |\nn P_{(t-\tau_n)^+}f|)(X_{t\land \tau_n})\\
 &\cdot \e^{-\ss 2 \int_0^{t \land\tau_n}\<u_s^{-1}\nn\log \phi (X_s),\, \d B_s\> -\int_0^{t\land\tau_n} (K+|\nn\log \phi|^2)(X_s)\d s}\Big\}.\end{split}\end{equation*}Since $K, |\nn\log\phi|$ are bounded and due to   Lemma \ref{L2.2}  $|\nn P_\cdot f|$ is bounded on
$[0,t]\times M$, according to the dominated convergence theorem we complete the proof by letting $n\to\infty$.
\end{proof}

\subsection{Equivalence of (2) and (3)} The equivalence of (2) and (3) follows from the following result by taking $\tt Z= -\ss 2\,\nn\log \phi$.
\beg{prp}\label{P2.3}  Let $\tt Z$ be a bounded $C^1$-vector field on $M$, and let $Y_t$ be the reflecting diffusion process generated by $L+\ss 2\, \tt Z$ starting at $x$. Then for any bound measurable function $F$ of $X_{[0,t]}:=\{X_s\}_{s\in [0,t]}$,
$$\E^x \Big\{F(X_{[0,t]})\e^{ \int_0^t\<u_s^{-1} \tt Z(X_s),\, \d B_s\>-\ff 1 2 \int_0^t |\tt Z|^2(X_s)\d s}\Big\} = \E^x F(Y_{[0,t]}).$$\end{prp} \beg{proof} Let
$$R= \exp\bigg[\int_0^t\<u_s^{-1} \tt Z(X_s),\, \d B_s\>-\ff 1 2 \int_0^t |\tt Z|^2(X_s)\d s\bigg].$$ By the Girsanov theorem, under the probability measure $R\d\P$ the process
$$\tt B_s:= B_s - \int_0^s \<u_r^{-1} \tt Z(X_r), \d B_r\>,\ \ s\in [0,t]$$ is a $d$-dimensional Brownian motion. Obviously, the equation (\ref{E}) can be reformulated as
$$\d X_s= \ss 2\, u_s\circ \d \tt B_s +\big(Z+\ss 2\, \tt Z\big)(X_s)\d s +N(X_s)\d l_s.$$ Therefore, under the new probability measure, $X_{[0,t]}$ is the reflecting diffusion process generated by $L+\ss 2\, \tt Z$. Hence,
$\E^x \{RF(X_{[0,t]})\} = \E^x F(Y_{[0,t]}).$\end{proof}

\subsection{From (3) to (1)} The desired assertion follows from the following result, which also indicates that
for a class of manifolds including compact ones, the condition $\II\ge -N\log\phi$ in the definition of $\D$ is essential for  (3), and hence (2).

\beg{prp}\label{P2.4} For any strictly positive function $\phi\in C_b^2(M)$, the gradient inequality in Theorem $\ref{T1.1}(3)$ implies $\Ric_Z^{\phi,1}\ge K$. If there exists $r_0>0$ such that on $\{\rr_\pp\le r_0\}$ the distance function $\rr_\pp$ to the boundary is smooth with bounded $L\rr_\pp$, then Theorem $\ref{T1.1}(3)$ also implies $\II\ge -N\log \phi$.\end{prp}

\beg{proof} (a) Let $x\in M\setminus \pp M$ and $X\in T_x M$ with $|X|=1$, we aim to prove $\Ric_Z^{\phi,1}(X,X)\ge K$ from   (3).
Let $f\in C_0^\infty(M)$ with supp$f\subset M\setminus \pp M$ be such that $\nn f(x)=X$ and $\Hess_f(x)=0.$ Let $\vv>0$ such that $|\nn f|\ge \ff 1 2$ on $B(x,\vv),$ the geodesic ball at $x$ with radius $\vv$. Let $X_t^\phi$ be the reflecting diffusion generated by $L^\phi$ with $X_0^\phi=x$, and let
$$\si_\vv=\inf\{t\ge 0: \rr(X_t^\phi,x)\ge \vv\}.$$ By Lemma\ref{L2.3},
$$\P(\si_\vv\ge t)\le \e^{-c/t},\ \ t\in (0,1]$$ holds for some constant $c>0$. Since  $l_s=0$ for $s\le \si_\vv$, this implies that
\beq\label{B2} \beg{split} &\E^x \Big\{(\phi|\nn f|)(X_t^\phi)\e^{\int_0^t K(X_s^\phi)\d s}\Big\}  = \E^x \Big\{(\phi|\nn f|)
(X^\phi_{t\land \si_\vv})\e^{\int_0^{t\land \si_\vv} K(X_s^\phi)\d s}\Big\} + {\rm o}(t)\\
&= (\phi|\nn f|)(x) + t \big\{L^\phi (\phi |\nn f|)+K \phi|\nn f|\big\}(x) +{\rm o}(t),\end{split}\end{equation} where ${\rm o}(t)$ stands for a $t$-dependent quantity such that $ {\rm o}(t)/t\to 0$ as $t\to 0$. On the other hand, since supp$f\subset M\setminus \pp M$ so that $Nf=0$, we have
$$P_tf= f+ \int_0^t P_s Lf\d s.$$ This and $\nn f(x)|=|X|=1$  imply that
$$|\nn P_t f(x)| =|\nn f +t(\nn Lf)|(x) +{\rm o}(t)= |\nn f|(x)+ \<\nn Lf,\nn f\>(x) t+ {\rm o}(t).$$ Combining this with (\ref{B2}) and the gradient inequality in Theorem \ref{T1.1}(3), we arrive at
\beq\label{FF}\big\{L^\phi(\phi|\nn f|)-\phi\<\nn L,\nn f\>\big\}(x)\ge -\big\{K\phi|\nn f|\big\}(x)=-(K\phi)(x).\end{equation}  Noting that $\Hess_f(x)=0$ and $|\nn f(x)|=1$ imply
$$L^\phi(\phi|\nn f|)(x) = L^\phi  \phi(x) +\{\phi L |\nn f|^2\}(x)=\{\phi(L|\nn f|^2 - \phi L\phi^{-1})\}(x),$$ combining (\ref{FF}) with (\ref{1.1}) we obtain $$\Ric_Z^{\phi,1}(X,X)= \Ric_Z^{\phi,1}(\nn f,\nn f)(x) \ge -K(x).$$

(b) Let $x\in \pp M$ and $X\in T_x\pp M$ with $|X|=1$. Let $f\in C_0^\infty(M)$ be such that $Nf=0$ and  $\nn f(x)=X$. We have
$$P_t f = f +\int_0^t P_s Lf \d s.$$ Consequently, for small $t$,
\beq\label{A2}  |\nn P_tf(x)|^2 = |\nn f(x)|^2+ {\rm o}(t^{1/2})= 1 +{\rm o}(t^{1/2}).\end{equation} On the other hand, according to
\cite[Proposition 4.1]{W10}, $$ \E^x l_t^\phi =\ff{2\ss t}{\ss\pi} + \circ (t^{1/2}),$$ where $l_t^\phi$ is the local time of $X_t^\phi$ on $\pp M$. Therefore,
since $|\nn f(x)|=1$ and
$$\lim_{s\to 0} N(\phi^2|\nn f|^2)(X_s^\phi) = N(\phi^2|\nn f|^2)(x),$$   we have \beg{equation*}\beg{split}
 &\E^x\Big\{ (\phi |\nn f|)(X_t^\phi) \e^{-\int_0^t  K (X_s^\phi)\d s}\Big\}^2=P_t^\phi(\phi |\nn f|)^2 (x) + {\rm o}(t^{1/2})\\
 &= (\phi^2|\nn f|^2)(x) +\int_0^t P_s^\phi L^\phi (\phi^2|\nn f|^2)(x)\d s +
\E^x\int_0^t N(\phi^2|\nn f|^2)(X_s^\phi)\d l_s^\phi\\
&= \phi^2(x) + \ff{2\ss t}{\ss\pi} N(\phi^2|\nn f|^2)(x) +{\rm o}(t^{1/2}).\end{split}\end{equation*}
Combining this with (\ref{A2}) and the gradient inequality in Theorem \ref{T1.1}(3), we conclude that
$$N(\phi^2|\nn f|^2)(x)\ge 0.$$ This implies  $\II(X,X)\ge - N\log \phi(x)$ since $X=\nn f(x)$ and by \cite[(3.8)]{W10b}, $N|\nn f|^2= 2 \II(\nn f,\nn f)$. \end{proof}

\section{Proofs of Corollaries \ref{C1.2} and \ref{C1.3}}

We first present two lemma    which are   known when $M$ is compact, where the first extends \cite[Theorem 2.1]{W09} and the second is crucial in order to use Bakry-Emery's semigroup argument. For readers' convenience we include below complete proofs for both of them. 

\beg{lem}\label{3.01} Let $x\in \pp M$ and let $\si_r$ be in Lemma $\ref{L2.3}$ for a fixed constant $r>0$. Then 
$$\limsup_{t\to 0} \ff 1 t \Big|\E l_{t\land\si_r}-\ff{2\ss{t}}{\ss\pi}\Big|<\infty.$$\end{lem} 

\beg{proof} The proof is modified from \cite{W09}. Let $\rr_\pp$ be the Riemannian distance to $\pp M$, and let $r_0\in (0, r)$ be such that $\rr_\pp$ is smooth on $B(x, 2 r_0).$ Let 
$$\tau=\inf\{t\ge 0: \rr(X_t, x)\ge r_0\}.$$ By the It\^o formula we have 
\beq\label{L1} \d \rr_\pp(X_t) =\ss 2\, \d b_t +L\rr_\oo(X_t)\d t+\d l_t,\ \ t\le\tau,\end{equation} where $b_t$ is an one-dimensional Brownian motion. Let $\tt B_t$ solve 
$$\d\tt b_t= \text{sgn}(\tt b_t)\d b_t,\ \ \tt b_0=0.$$ then $\tt b_t$ is an one-dimensional Brownian motion such that 
$$\d |\tt b_t| =\d b_t +\d\tt l_t,$$ where $\tt l_t$ is the local time of $\tt B_t$ at $0$. Combining this with (\ref{L1}) and noting that $\d l_t$ is supported on $\{\rr_\pp(X_t)=0\}$ while $\d\tt l_t$ is supported on $\{\tt b_t=0\}$, we obtain
\beg{equation*}\beg{split} &\d\big(\rr_\pp(X_t)-\ss 2\,|\tt b_t|\big)^2\\ &= 2\big(\rr_\pp(X_t)-\ss 2\,|\tt b_t|\big) L\rr_\pp (X_t)\d t+ 2 \big(\rr_\pp(X_t)-\ss 2\,|\tt b_t|\big)(\d l_t-\ss 2\,\d\tt l_t)\\
&\le 2\big(\rr_\pp(X_t)-\ss 2\,|\tt b_t|\big)L\rr_\pp(X_t)\d t\le c_1\big|\rr_\pp(X_t)-\ss 2\,|\tt b_t|\big|\d t,\ \ t\le\tau\end{split}\end{equation*} for some constant $c_1>0$. This implies 
$$\E\big(\rr_\pp(X_{t\land\tau})-\ss 2 \,|\tt b_{t\land\tau}|\big)^2 \le \ff {c_1^2} 4 t^2,\ \ t\ge 0.$$ Since due to 
(\ref{L1}) one has $|\E l_{t\land\tau} -\E \rr_\pp(X_{t\land\tau})\big|^2\le c_2t^2$ for some constant $c_2>0$, it follows that 
$$\big|\E l_{t\land\tau} -\ss 2\,\E|\tt b_{t\land \tau}|\big|\le c_3t,\ \ t\ge 0$$ holds for some constant $c_3>0.$ Noting that $\E|\tt b_t|=\ss{2t/\pi}$ and $\E \tt b_t^2=t$, combining this with Lemma \ref{L2.3} we arrive at 
\beq\label{L2} \beg{split}\Big|\E l_{t\land\tau}-\ff{2\ss t}{\ss\pi}\Big|&=\big|\E l_{t\land\tau}-\ss 2\,\E|\tt b_t|\big| \le c_3t +\ss 2 \E(|\tt b_t|1_{\{t>\tau\}})\\
&\le c_3 t +\ss{2t\P(t>\tau)}\le c_4 t,\ \ t\in [0,1]\end{split}\end{equation} for some constant $c_4>0.$ Finally, using $t\land \si_r$ in place of $t$, the proof of \cite[Lemma 2.3]{W09} leads to $\E l_{t\land \si_r}^2\le c_0 t$  for some constant $c_0>0$ and all $t\in [0,t].$ Therefore, it follows from (\ref{L2}) and Lemma \ref{L2.3} that
$$\Big|\E l_{t\land\si_r}-\ff{2\ss t}{\ss\pi}\Big|\le c_4 t \E(l_{t\land \si_r}1_{\{t>\tau\}})\le c_4t +\ss{c_0t \P(t>\tau)}\le c_5 t$$ holds for some constant $c_5>0$ and all $t\in [0,1].$ This completes the proof.\end{proof}

\beg{lem}\label{L3.0} Let  $f\in \C,$ the class of all functions $f\in C^2(M)$ such that $Nf=0$ and $Lf$ is bounded. Then
\beg{enumerate} \item[$(1)$] $\ff{\d}{\d t} P_t f= P_t Lf = LP_t f,\ t\ge 0;$
\item[$(2)$] $N P_t f=0,\ t\ge 0;$\item[$(3)$] Let $t>0$ and $F\in C_b^2([\inf f, \sup f])$. If $|\nn P_\cdot f|$ is bounded on $[0,t]\times M$, then $$\ff{\d}{\d s} P_s F(P_{t-s}f)= P_s \big(F''(P_{t-s}f)|\nn P_{t-s}f|^2\big),\ \ s\in [0,t].$$\end{enumerate}\end{lem}
\beg{proof} (1) The first equality follows from $P_t f= f +\int_0^t P_s Lf\d s$ implied by the It\^o formula. To prove the second equality, it suffices to show that for any $x\in M\setminus \pp M$,
\beq\label{3.01} \ff{\d}{\d t} P_t f(x)= LP_t f(x).\end{equation} Let $r_0>0$ be such that $B(x, r_0)\subset M\setminus \pp M$, and take $h\in C_0^\infty(M)$ such that $h|_{B(x, r_0/2)}=1$ and $h|_{B(x,r_0)^c}=0.$ By the It\^o formula we have
$$P_{t+s}f(x)= \E^x (h P_t f)(X_s) + \E^x \{(1-h) P_t f\}(x) = \E^x\int_0^s L(h P_t f)(X_r)\d r + \E^x \{(1-h) P_t f\}(x).$$ Since $L(h P_t f)(X_r)$ is bounded and goes to $LP_t f(x)$ as $r\to 0$, and noting that by Lemma \ref{L2.3}
$$\E^x |(1-h)P_t f|(X_s)\le \|f\|_\infty \e^{-c/s},\ \ s\in (0,1]$$ holds for some constant $c>0$, we conclude that
$$\ff{\d }{\d t} P_t f(x)= \lim_{s\downarrow 0} \ff{P_{t+s}f(x)-P_t f(x)} s = LP_t f(x),$$ that is, (\ref{3.01}) holds.

(2) Let $x\in \pp M$. If $N P_t f(x) \ne 0$, for instance $N P_t f(x)>0$, then there exists $r_0,\vv>0$ such that $NP_t f\ge \vv$ holds on $B(x, 2 r_0)$. Moreover, by using $f+\|f\|_\infty$ in place of $f$, we may assume that $f\ge 0$.  Let $h\in C_0^\infty(M)$ such that $0\le h\le 1, N h=0, h|_{B(x,r_0)}=1$ and $h|_{B(x, 2 r_0)^c} =0$. By the It\^o formula and using (1), we obtain
\beg{equation*}\beg{split} P_{t+s} f(x)&\ge P_s (hP_t f)(x) = P_t f(x) +\int_0^s P_r L(h P_t f)(x)\d r +\E^x \int_0^s (h NP_t f)(X_r) \d r\\
 &\ge P_t f(x) + s P_{t} Lf(x) +\text{o}(s) + \vv \E^x l_{s\land\si},\end{split}\end{equation*} where $\si:=\inf\{s\ge 0: X_s\notin B(x, r_0)\}$. Combining this with (1) we arrive at
\beq\label{CC}\vv \lim_{s\to 0}\ff 1 s \E^x l_{s\land\si}\le 0,\end{equation} which is impossible  according to Lemma \ref{L2.3} and \cite[Theorem 2.1]{W09}. Indeed, since
$l_{s\land \tau}$ only depends on the process before exiting $B(x,r_0)$, we may assume that $M$   is compact (otherwise, simply use a large enough smooth bounded domain   to replace $M$). In this case \cite[Theorem 2.1]{W09} implies that $E^x l_s\ge c\ss s$ holds for some constant $c>0$ and small $s>0$. Thus, it follows from Lemma \ref{L3.0} and \cite[Lemma 2.2]{W09} that
$$\E^x l_{s\land\si}\ge \E^x l_s - (\E^x l_s^2)^{1/2} \P(\si<s)^{1/2} \ge c'\ss s$$ holds for some constant $c'>0$ and small $s>0$. This is contradictive to  (\ref{CC}).

(3) By (1) and (2) and using the It\^o formula, there is a local martingale $M_s$ such that
\beg{equation*}\beg{split} \d F(P_{t-s}f)(X_s) &= \d M_s \{ LF(P_{t-s} f)+F'(P_{t-s}f)LP_{t-s}f\}(X_s)\d s \\
&= \d M_s + \{F''(P_{t-s} f)  |\nn P_{t-s}f|^2\}(X_s)\d s,\ \ s\in [0,t].\end{split}\end{equation*} Since $|\nn P_\cdot f|$ is bounded on $[0,t]\times M$ and $F\in C_b^2([\inf f, \sup f])$, we see that $M_s$ is indeed a martingale. Therefore,
$$P_s F(P_{t-s}F)= \E F(P_{t-s}f)(X_s)= F(P_t f) +\int_0^s P_r \{F''(P_{t-s}f)|\nn P_{t-s}f|^2\}\d r.$$ This completes the proof.
\end{proof}

Next, we present a result on the   Poincar\'e type inequalities and the log-Harnack inequality for $P_t$ by using an $L^2$-gradient estimate. Having Lemma \ref{L3.0} in hands,    the proof of (\ref{P1}) is standard according to Bakry and Ledoux (see e.g. \cite{Bakry, Ledoux}), while that of (\ref{LH1}) is essentially due to \cite{RW}.

\beg{lem}\label{L3.1} If $|\nn P_t f|^2 \le \xi_t P_t |\nn f|^2$ holds for some strictly positive $\xi\in C([0,\infty))$ and all $t\ge 0$ and $f\in C_b^1(M)$, then
\beq\label{P1} 2|\nn P_t f|^2 \int_0^t\ff  {\d s} {\xi_s}\le  P_t f^2-(P_tf)^2\le 2(P_t|\nn f|^2)\int_0^t \xi_s\d s,\ \ t\ge 0, f\in C_b^1(M),\end{equation} and for any measurable function $f$ with $f\ge 1$, \beq\label{LH1} P_t\log f(y)\le \log P_t f(x) +\ff{\rr(x,y)^2}{4\int_0^t\xi_s^{-1}\d s},\ \ t>0.\end{equation} \end{lem}

\beg{proof} It suffices to prove for $f\in \C$ such that $f\ge 1$. For any $\vv>0$ let   $\gg: [0,1]\to M$ be the minimal curve such that $\gg(0)=x, \gg(1)=y$ and $|\dot\gg|\le \rr(x,y)+\vv.$ Let
$$h(s)= \ff{\int_0^s \xi_r^{-1}}{\int_0^t \xi_r^{-1}\d r},\ \ s\in [0,t].$$ By Lemma \ref{L3.0},
we have
\beg{equation*}\beg{split} &\ff{\d}{\d s} (P_s\log P_{t-s} f)(\gg\circ h(s))\\
 &\le -P_s |\nn\log P_{t-s} f|^2(\gg\circ h(s)) +(\vv+\rr(x,y))\dot h(s) |\nn P_s\log P_{t-s}f|(\gg\circ h(s))\\
&\le\big\{ -P_s |\nn\log P_{t-s} f|^2 +(\vv+\rr(x,y))\dot h(s) \ss{\xi_s P_s |\log P_{t-s}f|^2}\big\}(\gg\circ h(s))\\
&\le \ff 1 4 \rr(x,y)^2 \dot h(s)^2 = \ff{(\vv+\rr(x,y))^2}{4\xi_s(\int_0^t \xi_r^{-1}\d r)^2},\ \ \ s\in (0,t).\end{split}\end{equation*} Integrating over $[0,t]$ and letting $\vv\downarrow 0$, we obtain (\ref{LH1}).

Next, noting that
$$\ff{\d}{\d s} P_s (P_{t-s}f)^2 = 2 P_s |\nn P_{t-s} f|^2\beg{cases}\le 2\xi_{t-s} P_t|\nn f|^2,\\
\ge \ff{2}{\xi_s}|\nn P_t f|^2,\end{cases} \ \ s\in (0,t),$$  we prove (\ref{P1}). \end{proof}

\beg{proof}[Proof of Corollary \ref{C1.2}] Due to Lemma \ref{L3.1} and $\phi\ge 1$, it suffices to prove the first assertion. Obviously, $\Ric_Z^{\phi,2}\ge K_\phi$ implies that $\Ric_Z^{\phi,1}\ge K:= K_\phi +|\nn\log \phi|^2.$ Let
$$R_t= \exp\bigg[ -\ss 2\int_0^t \<u_s^{-1}\nn\log \phi(X_s), \d B_s\> -\int_0^t |\nn\log \phi(X_s)|^2\d s\bigg].$$  By Theorem \ref{T1.1} and $\phi\ge 1$, we obtain
\beg{equation*}\beg{split} (\phi |\nn P_t f |)^2(x) &\le \Big(\E^x \Big\{R_t (|\nn f|\phi)(X_t)\e^{-\int_0^t K(X_s)\d s}\Big\}\Big)^2\\
&\le \{P_t(\phi |\nn f|)^2(x)\} \E^x \Big(R_t^2 \e^{-2\int_0^t K(X_s)\d s}\Big)\\
&\le\{P_t(\phi |\nn f|)^2(x)\} \e^{-2K_\phi t} \E^x \e^{-2\ss 2 \int_0^t \<u_s^{-1}\nn\log \phi(X_s),\, \d B_s\> -4\int_0^t |\nn\log \phi(X_s)|^2\d s }\\
& =  \e^{-2K_\phi t} P_t(\phi |\nn f|)^2(x).\end{split}\end{equation*}\end{proof}

To prove Corollary \ref{C1.3}, we present a  log-Sobolev inequality which generalizes the corresponding known one on manifolds without boundary.

\beg{lem}\label{L3.2} Let $\phi\in\D$ such that $\Ric_Z^{\phi,2}\ge K_\phi$ holds for some constant $K_\phi$. Let  $\bar P_t^\phi$ be  the  semigroup of the reflecting diffusion process generated by $\bar L^\phi:=L-4\nn\log\phi$. Then
\beq\label{LS1} P_t(f^2\log f^2) \le (P_t f^2)\log P_t f^2 + 4 \|\phi\|_\infty^2 \int_0^t \e^{-2K_\phi(t-s)} P_s \bar P_{t-s}^\phi|\nn f|^2\ d s\end{equation} holds for all $t\ge 0$ and $f\in C_b^1(M).$\end{lem}

\beg{proof} It suffices to prove for $f\in\C$ with $\inf f^2>0$. Let $R_t$ be in the proof of Corollary \ref{C1.2}.
Since $\Ric_Z^{\phi,2}\ge K_\phi$ implies that $\Ric_Z^{\phi,1}\ge K:= K_\phi +|\nn\log\phi|^2$, by Theorem \ref{T1.1} and $\phi\ge 1$ we have
\beg{equation*}\beg{split} |\nn P_t f^2(x)|^2 &\le \Big(\E^x\Big\{R_t (\phi |\nn f^2|)(X_t)\e^{-\int_0^t K(X_s)\d s}\Big\}\Big)^2\\
&\le 4\|\phi\|_\infty^2 (P_t f^2(x)) \E^x\Big\{R_t^2|\nn f|^2(X_t)\e^{-2\int_0^t K(X_s)\d s}\Big\}\\
&\le 4\|\phi\|_\infty^2 (P_tf^2(x))\E^x\Big\{|\nn f|^2(X_t) \e^{-2\ss 2  \int_0^t \<u_s^{-1}\nn\log \phi(X_s),\, \d B_s\> -2\int_0^t (|\nn\log \phi(X_s)|^2+K(X_s))\d s }\Big\}\\
&= 4\|\phi\|_\infty^2 (P_t f^2(x))\e^{-2K_\phi t} \E^x  \{\bar R_t |\nn f|^2(X_t)\},\end{split}\end{equation*} where
$$\bar R_t:= \e^{-2\ss 2  \int_0^t \<u_s^{-1}\nn\log \phi(X_s),\, \d B_s\> -4\int_0^t |\nn\log \phi(X_s)|^2 \d s }.$$
 Combining this with Proposition \ref{P2.3} for $\tt Z= -2\ss 2\, \nn\log\phi$, we obtain
$$|\nn P_t f^2|^2 \le 4 \|\phi\|_\infty^2 (P_tf^2)\e^{-2K_\phi t} \bar P_t^\phi|\nn f|^2,\ \ t\ge 0.$$
 Therefore, by Lemma \ref{L3.0},
$$\ff{\d}{\d s} P_s\big\{(P_{t-s}f^2)\log P_{t-s}f^2\big\}=P_s\ff{|\nn P_{t-s}f^2|^2}{P_{t-s} f^2 }\le 4\|\phi\|_\infty^2 \e^{-2K_\phi(t-s)}P_s\bar P_{t-s}^\phi|\nn f|^2.$$ Then the proof is  completed  by integrating   over $[0,t]$.\end{proof}

\beg{proof}[Proof of Corollary \ref{C1.3}]  Let $f\in C_b^1(M)$ such that $\mu(f^2)=1$ and $\mu(|\nn f|^2)>0$.
Since  $\mu$ is $P_t$-invariant while $\phi^{-4}\d\mu$ is   $\bar P_t^\phi$-invariant, integrating (\ref{LS1}) w.r.t. $\mu$ gives
\beq\label{LS2}\beg{split} \mu(f^2\log f^2) &\le \mu((P_tf^2)\log P_t f^2) + 4 \|\phi\|_\infty^2 \int_0^t\e^{-2K_\phi s}\mu(\bar P_s^\phi |\nn f|^2)\d s\\
&\le   \mu((P_tf^2)\log P_t f^2) + 4 \|\phi\|_\infty^6 \int_0^t\e^{-2K_\phi s}\mu(\phi^{-4} |\nn f|^2)\d s\\
&\le   \mu((P_tf^2)\log P_t f^2) + \ff{2\|\phi\|_\infty^6(1-\e^{-2K_\phi t})}{K_\phi} \mu( |\nn f|^2).\end{split}\end{equation} If $K_\phi>0$, then letting $t\to\infty$ we prove Corollary \ref{C1.3}(1).

Now, the proof of the second assertion can be done as in \cite{BGL} by using (\ref{LS2}) and Corollary \ref{C1.2}(2).  Applying Corollary \ref{C1.2}(2) for $P_t f^2$ in place of $f$, we find
$$P_t \log P_t f^2(y)\le \log P_{2t}f^2(x)+\ff{\|\phi\|_\infty^2K_\phi \rr(x,y)^2}{2(\e^{2K_\phi t}-1)},\ \ x,y\in M, t> 0.$$ Integrating w.r.t. the optimal coupling of $f^2\mu$ and $\mu$, which reaches the inf in the definition of $W_2^\rr(f^2\mu,\mu)$, and noting that $P_t$ is symmetric in $L^2(\mu)$,  we obtain
$$\mu((P_t f^2)\log P_t f^2) \le \ff{\|\phi\|_\infty^2K_\phi W_2^\rr(f^2\mu,\mu)^2}{2(\e^{2K_\phi t}-1)}.$$ Combining this with the first inequality in (\ref{LS2}), we arrive at
\beq\label{FW} \mu(f^2\log f^2) \le   \|\phi\|_\infty^6\mu(|\nn f|^2) r_t +\ff {\|\phi\|_\infty^2} {r_t}  W_2^\rr(f^2\mu,\mu)^2 -\ff{\|\phi\|_\infty^2K_\phi} 2 W_2^\rr(f^2\mu,\mu)^2,\end{equation} where $$r_t:= \ff{2(1-\e^{-2K_\phi t})} {K_\phi},\ \ t\ge 0.$$ If $K_\phi\le 0$ then $\{r_t:\ t\in [0,\infty]\}=[0,\infty]$.   So, there exists  $t\in [0,\infty]$ such that $$r_t =\ff{W_2^\rr(f^2\mu,\mu)}{\|\phi\|_\infty^2\ss{\mu(|\nn f|^2)}}.$$ Therefore, the desired HWI inequality follows from (\ref{FW}).
 \end{proof}

\beg{thebibliography}{99}

\bibitem{ATW} M. Arnaudon,  A. Thalmaier,    F.-Y. Wang,
 \emph{ Gradient estimates and Harnack inequalities on non-compact Riemannian manifolds,} Stoch. Proc. Appl. 119(2009), 3653--3670.

\bibitem{Bakry} D. Bakry, \emph{On Sobolev and logarithmic Sobolev
inequalities for Markov semigroups,} ``New Trends in Stochastic
Analysis'' (Editors: K. D. Elworthy, S. Kusuoka, I. Shigekawa),
Singapore: World Scientific, 1997.

 \bibitem{BE}   D. Bakry,  M. Emery, \emph{Hypercontractivit\'e de
semi-groupes de diffusion}, C. R. Acad. Sci. Paris. S\'er. I Math.
299(1984), 775--778.

\bibitem{Bess} A. L. Bess, \emph{Einstein Manifolds,} Sringer, Berlin, 1987.

\bibitem{BGL} S. G. Bobkov, I.  Gentil and M.  Ledoux,
\emph{Hypercontractivity of Hamilton-Jacobi equations,} J. Math.
Pures Appl. 80(2001), 669--696.

 \bibitem{Hsu} E. P. Hsu, \emph{Multiplicative functional for the heat
equation on manifolds with boundary,} Michigan Math. J. 50(2002),
351--367.

\bibitem{Ledoux} M. Ledoux, \emph{The geometry of Markov diffusion generators,} Ann. Facu. Sci. Toulouse
9(2000), 305--366.

\bibitem{OV} F.  Otto,  C. Villani, \emph{Generalization of an inequality
by Talagrand and links with the logarithmic Sobolev inequality,} J.
Funct. Anal.  173(2000), 361--400.

\bibitem{RW}  M. R\"ockner, F.-Y. Wang, \emph{Log-Harnack  inequality for stochastic differential equations in Hilbert spaces and its consequences, } Infin. Dimens. Anal. Quant. Probab.  Relat. Topics 13(2010), 27--37.


 \bibitem{W05} F.-Y. Wang, \emph{Gradient estimates and the first
Neumann eigenvalue on manifolds with boundary,} Stoch. Proc. Appl.
115(2005), 1475--1486.

\bibitem{W07} F.-Y. Wang, \emph{Estimates of the first Neumann eigenvalue and the log-Sobolev constant
  on  nonconvex manifolds,} Math. Nachr. 280(2007), 1431--1439.
\bibitem{W09}  	F.-Y. Wang, \emph{Second fundamental form and gradient of Neumann semigroups,}  J. Funct. Anal. 256(2009), 3461--3469.
\bibitem{W09b} F.-Y. Wang, \emph{Log-Sobolev inequalities: different roles of Ric and Hess,}  Annals of Probability 37(2009), 1587--1604.
\bibitem{W09c} 	F.-Y. Wang, \emph{Log-Sobolev inequality on non-convex manifolds,}   Adv. Math. 222(2009), 1503-1520.

\bibitem{W10}	F.-Y. Wang, \emph{Harnack inequalities on manifolds with boundary and applications,}  J. Math. Pures Appl. 94(2010), 304--321.

\bibitem{W10b} F.-Y. Wang, \emph{Semigroup properties for the second fundamental form,} Docum. Math. 15(2010), 527--543.
\end{thebibliography}

\end{document}